\documentclass{article}[11pt]
\usepackage{amssymb}
\usepackage{amsmath}
\DeclareMathOperator{\ch}{ch}
\DeclareMathOperator{\sh}{sh}
\DeclareMathOperator{\tr}{tr}
\newcommand{\Cov}{\mathrm{Cov}}
\newcommand{\der}{\mathrm{d}}
\newcommand{\dvol}{\mathrm{d}\,\mathrm{vol}}

\newcommand{\Img}{\mathrm{Im}}
\newcommand{\Ker}{\mathrm{Ker}}
\newcommand{\M}{\mathcal{M}}
\newcommand{\Mn}{\mathbb{E}}
\newcommand{\N}{\mathbb{N}}
\newcommand{\nep}{\mathrm{e}}
\newcommand{\Pb}{\mathbb{P}}
\newcommand{\Ric}{\underline{\mathrm{Ric}}}
\newcommand{\rk}{\mathrm{rk}}
\newcommand{\Rl}{\mathbb{R}}
\newcommand{\T}{\mathrm{T}}
\newcommand{\Var}{\mathrm{Var}}
\newtheorem{thm}{Theorem}
\newtheorem{rmk}[thm]{Remark}
\newtheorem{lem}[thm]{Lemma}
\newtheorem{xpl}[thm]{Example}
\newtheorem{defi}[thm]{Definition}
\title{Improved spectral gap bounds on positively curved manifolds}
\author{Laurent Veysseire}
\begin{document}
\maketitle
\abstract{A coupling method and an analytic one allow us to prove new lower bounds for the spectral gap of reversible diffusions on compact manifolds.
Those bounds are based on the a notion of curvature of the diffusion, like the coarse Ricci curvature or the Bakry--Emery curvature-dimension inequalities.
We show that when this curvature is nonnegative, its harmonic mean is a lower bound for the spectral gap.}
\section*{Introduction}

The study of the spectrum of the Laplace Operator on Riemannian manifolds has many applications in various domains of mathematics.
A whole chapter of \cite{berg} is devoted to this issue.
In this article, we take the convention
$$\Delta =g^{ij}\nabla_i\nabla_j$$
for the Laplace operator.
The spectral gap of $\Delta$ is the opposite of the greatest non-zero eigenvalue of $\Delta$ (the spectrum of $\Delta$ is discrete and non-positive).

One way to estimate this spectral gap is to use the Ricci curvature, as we see it in the Lichnerowicz theorem (see \cite{lich}).
\begin{thm}[Lichnerowicz]\label{Lic}
Let $\M$ be a $n$-dimensional Riemannian manifold.
If there exists $K>0$ such that for each $x\in\M$, for each $u\in\T_x\M$, we have $\mathrm{Ric}_x(u,u)\geq K g_x(u,u)$, then the spectral gap $\lambda_1$ of $\Delta$ satisfies
$$\lambda_1\geq\frac{n}{n-1}K.$$
Here we denote by $\mathrm{Ric}$ the Ricci curvature of $\M$.
\end{thm}

Chen and Wang improved this result in \cite{chw}, using the diameter of the manifold in their estimates:
\begin{thm}Let $\M$ be a compact connected $n$-dimensional Riemannian manifold, $K$ be the infimum of the Ricci curvature on $\M$ and $D$ be the diameter of $\M$.
Then if $K\geq0$, we have the following bounds:
$$\lambda_1\geq\frac{\pi^2}{D^2}+\max\left(\frac{\pi}{4n},1-\frac{2}{\pi}\right)$$
and if $n>1$,
$$\lambda_1\geq\frac{nK}{(n-1)\left(1-\cos^n\left(\frac{D\sqrt{K(n-1)}}{2}\right)\right)}.$$
And if $K\leq0$, we have the following bounds:
$$\lambda_1\geq\frac{\pi^2}{D^2}+\left(\frac{\pi}{2}-1\right)K$$
and if $n>1$,
$$\lambda_1\geq\frac{\pi^2\sqrt{1-\frac{2D^2K}{\pi^4}}}{D^2\ch\left(\frac{D\sqrt{-K(n-1)}}{2}\right)}$$
\end{thm}

In \cite{aubr}, E.Aubry gives a lower bound for $\lambda_1$ when the curvature
$$\Ric(x):=\inf_{u\in\T_x\M}\frac{\mathrm{Ric}_x(u,u)}{\|u\|^2}$$
is close to a positive constant in the sense of $L^p$ norm with $p$ large enough:

\begin{thm}Let $\M$ be a complete $n$-dimensional Riemannian Manifold, $p>\frac{n}{2}$ and $K>0$, such that
$$\int_{\M}(\Ric-K)_-^p<+\infty.$$
Then $\M$ has a finite volume and the spectral gap of $\Delta$ on $\M$ satisfies:
$$\lambda_1\geq\frac{n}{n-1}K\left(1-\frac{C(p,n)}{K}\|(\Ric-K)_-\|_p\right)$$
where $C(p,n)$ is a constant only depending on $p$ and $n$, and $\|f\|_p=\left(\frac{\int_{\M}|f|^p}{\mathrm{vol}(\M)}\right)^{\frac{1}{p}}$.
\end{thm}

This allows a little negative curvature, which is not the case of our results.

This article recapitulates and extends the results already stated in \cite{vey} and presents a coupling method, more adapted to discrete spaces than the analytic one.

We show by a coupling method that another bound for $\lambda_1$ is the harmonic mean of the Ricci curvature.
\begin{thm}\label{bharm} Let $\M$ be a compact Riemannian manifold with positive Ricci curvature.
Then we have $$\frac{1}{\lambda_1}\leq\int_\M\frac{\der\mu(x)}{\Ric(x)},$$
with $\der\mu=\frac{\dvol}{\mathrm{vol}(\M)}$, where $\mathrm{vol}$ is the Riemannian volume measure on $\M$.\end{thm}
This bound is often better than the Lichnerowicz one because the harmonic mean is better (and can be much better) than the infimum.
But unfortunately we lose the $\frac{n}{n-1}$ factor.

Merging the proof of Theorem \ref{Lic} and an analytic proof of Theorem \ref{bharm} gives us the following improvement:
\begin{thm}\label{harm2} Let $\M$ be a Riemannian manifold with positive Ricci curvature and $K=\inf_{x\in\M}\Ric(x)$.
Then for every $0\leq c\leq K$, we have:
$$\lambda_1\geq\frac{n}{n-1}c+\frac{1}{\int_\M\frac{\dvol}{\Ric(x)-c}}.$$
\end{thm}

Taking $c=K$ gives us the Lichnerowicz bound or even better, while $c=0$ gives us Theorem \ref{bharm}.

Our coupling approach is based on a notion of coarse Ricci curvature, introduced by Yann Ollivier in \cite{olli}, which uses the Wasserstein distance $W_1$.
A major step in our proof is the use of the coupling given by the following theorem:
\begin{thm}\label{cent} Let $\M$ be a smooth Riemannian manifold, and $F^i$ be a smooth vector field on $\M$.
Assume that there exists a diffusion process associated with the generator $Lf=\Delta f+F^i\nabla_if$.
Let $\kappa(x,y)$ be the coarse Ricci curvature of the diffusion between $x$ and $y$ (see Definition \ref{kxy}).
Then for any two distinct points $x$ and $y$ of $\M$, there exists a coupling $(x(t),y(t))$ between the paths of the diffusion process starting at $x$ and $y$ which satisfies:
$$d(x(t),y(t))=d(x,y)\nep^{-\int_0^t\kappa(x(s),y(s))\der s}$$
on the event that for any $s\in [0,t]$, $(x(s),y(s))$ does not belong to the cut-locus of $\M$.
\end{thm}

The contraction rate $\kappa(x,y)$ of this coupling behaves like the one of the coupling derived from the diffusion in $\mathcal{C}^1$ path space defined by M.Arnaudon, K.A.Coulibaly and A.Thalmaier in \cite{arn} when $x$ and $y$ are close.
We have a cut-locus problem that we will avoid by making a compactness assumption, which was anyway necessary to replace $\kappa(x,y)$ by its limit when $x$ and $y$ are infinitely close.

The coupling method and the analytic one keep working when we add a drift to the Brownian motion, provided the diffusion is reversible.
In this case, the generator takes the following form:
$$L=\frac{1}{2}g^{ij}(\nabla_i\nabla_j-(\nabla_j\varphi)\nabla_i)$$
with $\varphi$ a smooth function on $\M$, and $\nep^{-\varphi}\dvol$ is a reversible measure.
We have then the following generalization of Theorem \ref{harm2}:

\begin{thm}\label{harm3}Let $\M$ be a compact Riemannian manifold and $L=\frac{1}{2}g^{ij}(\nabla_i\nabla_j-(\nabla_j\varphi)\nabla_i)$ be the operator associated with a reversible diffusion process on $\M$.
Suppose that we have a curvature-dimension inequality in the sense of Bakry-\'Emery (see \cite{bakem} or \cite{bakem2}) with a positive curvature $\rho$ and a constant and positive dimension $n'$,which is
$$\Gamma_2(f)(x)\geq\rho(x)\Gamma(f)(x)+\frac{1}{n'}L(f)(x)^2.$$
Let $R$ be the infimum of $\rho$.
Then for every $0\leq c<R,$ we have
$$\lambda_1(L)\geq\frac{n'}{n'-1}c+\frac{1}{\int_\M\frac{\der\pi(x)}{\rho(x)-c}}$$
with $\der\pi=\frac{\nep^{-\varphi}\dvol}{\int_\M\nep^{-\varphi}\dvol}$ the reversible probability measure.
\end{thm}

We try to generalize our coupling method to diffusions which are not adaptated to the metric $g$, that is, whose generator takes the more general form:
$$L=\frac{1}{2}A^{ij}\nabla_i\nabla_j+F^i\nabla_i$$
without having necessarily $A^{ij}=g^{ij}$ anymore.
We have a generalization of Theorem \ref{cent} only on the very restrictive condition:
$$(H)\Leftrightarrow\forall u\in\T\M,u^ig_{jk}u^jg_{lm}u^l\nabla_iA^{km}=0\Leftrightarrow g^{il}\nabla_lA^{jk}+g^{jl}\nabla_lA^{ki}+g^{kl}\nabla_lA^{ij}=0$$
and with a lower $\tilde{\kappa}(x,y)$ instead of $\kappa(x,y)$.
Note that $(H)$ is true for $A^{ij}=g^{ij}$, in which case we have $\tilde{\kappa}=\kappa$.

We have the following generalization of Theorem \ref{bharm}:
\begin{thm}\label{harm} Consider a diffusion process on a compact Riemannian manifold $\M$ which is reversible and satisfies $(H)$.
For every $x$ in $\M$, we set $\tilde{\kappa}(x)=\inf_{u\in\T_x\M}\tilde{\kappa}(x,u)$.
If we have $\tilde{\kappa}(x)\geq\varepsilon>0$, then the spectral gap of $L$ is at least the harmonic mean of $\tilde{\kappa}$ (with respect to the reversible probability measure $\pi$):
$$\frac{1}{\lambda_1(L)}\leq\int_\M\frac{\der\pi(x)}{\tilde{\kappa}(x)}.$$
\end{thm}

In section \ref{arg}, we present a short argument which shows how we can derive the harmonic mean from Theorem \ref{cent}.
In section \ref{rcm}, we define the coarse Ricci curvature for diffusions and construct our couplings, so it's where Theorem \ref{cent} is proved.
In section \ref{bnd}, we present the proofs using the couplings and purely analytical ones for the harmonic mean bounds for the spectral gap.

\section{The harmonic mean in a nutshell}\label{arg}

The result and its proof presented in this section are a shortcut found by Yann Ollivier to obtain a harmonic mean from Theorem \ref{cent}.

Using a classical method, we will prove thanks to Theorem \ref{cent} the following result, which is a weaker version of Theorem \ref{bharm}:
\begin{thm}\label{sim} Let $\M$ be a compact Riemannian manifold with positive Ricci curvature, and $f$ be any $1$-Lipschitz function on $\M$.
Then the variance of $f$ is at most the average of $\frac{1}{\Ric}$.
\end{thm}

Indeed, the Poincar\'e inequality states that $\Var_\mu(f)\leq\frac{1}{\lambda_1}\int \|\nabla f\|^2\der\mu$, and the integral on the right hand side is at most $1$ for $1$-Lipschitz functions. In \cite{mil}, E.Milman shows that the converse is true i.e a control on the variance of Lipschitz functions (and even on the $L^1$ norm of $0$-mean Lipschitz functions) implies a Poincar\'e inequality, with a universal loss in the constants, under the hypothesis of a Bakry--Emery $CD(0,\infty)$ curvature-dimension inequality.

\noindent\textbf{Proof of theorem \ref{sim}: }
We only have to prove the result for $f$ regular enough, and use a density argument to get the result for non-regular $f$.
We consider the semi-group $P^t$ generated by the Laplacian operator. Then the limit of $P^t$ when $t$ tends to infinity is the operator which associates to $f$ the constant function equals to the mean of $f$ (respect to the normalized Riemannian volume measure).
So the variance of $f$ is the limit of the mean of $P^t(f^2)-(P^t(f))^2$ when $t$ tends to infinity.
We have
$$P^t(f^2)-(P^t(f))^2=\int_0^t\frac{\der}{\der s}\left(P^s((P^{t-s}(f))^2)\right)\der s=\int_0^tP^s(2\|\nabla(P^{t-s}(f))\|^2)\,\der s$$
Integrating over $\M$ yields
$$\int_\M(P^t(f^2)(x)-(P^t(f)(x)))^2\,\dvol(x)=\int_\M\int_0^t2\|\nabla(P^{t-s}(f))(x)\|^2\,\der s\,\dvol(x).$$
Thanks to Theorem \ref{cent}, by taking $y$ very close to $x$, we have $\|\nabla(P^{t-s}(f))(x)\|\leq \Mn_{\Pb_x}[\nep^{-\int_0^{t-s}\Ric(X_u)\der u}\|\nabla f(X_{t-s})\|]$, where the right hand side is the expectation of the term inside the brackets when $X$ has the law $\Pb_x$ of the twice accelerated Brownian motion on $\M$ starting at $x$.
Using the convexity of the exponential function, and the fact that $f$ is $1$-Lipschitz, we get then
$$\begin{array}{r@{\,}c@{\,}l}\int_\M(P^t(f^2)(x)-(P^t(f)(x)))\dvol(x)&\leq&2\int_\M\int_0^t\left(\Mn_{\Pb_x}\left[\nep^{-\int_0^{t-s}\Ric(X_u)\der u}\right]\right)^2\der s\,\dvol(x)\\
{}&\leq&2\int_\M\int_0^t\Mn_{\Pb_x}\left[\nep^{-2\int_0^{t-s}\Ric(X_u)\der u}\right]\der s\,\dvol(x)\\
{}&\leq&2\int_0^t\int_\M\Mn_{\Pb_x}\left[\int_0^1\nep^{-2(t-s)\Ric(X_{(t-s)u})}\der u\right]\dvol(x)\,\der s\\
{}&=&2\int_0^t\int_\M\nep^{-2(t-s)\Ric(x)}\dvol(x)\,\der s\\
{}&=&\int_\M\frac{1-\nep^{-2t\Ric(x)}}{\Ric(x)}\dvol(x).\end{array}$$
We just have to take the limit when $t$ tends to infinity and to divide by $\int_\M\dvol(x)$ to get the theorem.$\square$

\section{Coarse Ricci curvature for diffusions on Riemannian manifolds}\label{rcm}
In this section, we introduce the Coarse Ricci curvature $\kappa$ for general diffusions and give an explicit formula.
Then we construct the coupling of Theorem \ref{cent}, we show why the $(H)$ condition is needed and we define $\tilde{\kappa}$ when it is satisfied.

\subsection{Coarse Ricci curvature: definition and calculation}

Following what is done in \cite{olli} for Markov chains, we define the coarse Ricci curvature of diffusions as the rate of decay of the Wasserstein distance $W_1$ between the measures associated with the diffusion and starting at two different points:

\begin{defi}\label{kxy} Let $\M$ be a Riemannian manifold and $P^t$ be the semi-group of a diffusion on $\M$.
The coarse Ricci curvature between two different points $x$ and $y$ is the following quantity:
$$\kappa(x,y)=\varliminf_{t\rightarrow0}\frac{d(x,y)-W_1(\delta_x.P^t,\delta_y.P^t)}{td(x,y)}.$$
\end{defi}

The Wasserstein distance $W_1$ between two measures is the infimum on all the couplings of the expectation of the distance.
Our coupling will be consrtucted thanks to optimal ones.

To get an expression of this curvature only depending on the coefficients of the generator of the diffusion, we need to make sure that the diffusion does not move far away too fast.

\begin{defi} A diffusion on $\M$ is said to be locally uniformly $L^1$-bounded at $x$ if $\exists M>0,\exists\eta>0,\forall y\in\M|d(x,y)<\eta,\forall 0<t<\eta,\int d(x,z)\der (\delta_y.P^t)(z)<M$.
\end{defi}

\begin{rmk}If $\M$ is compact, any diffusion is locally uniformly $L^1$-bounded at each point (it suffices to take $M$ equals to the diameter of $\M$ in the previous definition).
\end{rmk}

The following theorem gives an expression of $\kappa(x,y)$.
Recall that the generator of the diffusion is
$$L(f)=\frac{1}{2}A^{ij}\nabla_i\nabla_jf+F^i\nabla_if$$
with $A$ symmetric and non-negative.
\begin{thm}\label{rxy} Take two distinct points $x$ and $y$ in $\M$, such that $A$ and $F$ are continuous at $x$ and $y$, and that the diffusion is locally uniformly $L^1$-bounded at $x$ and $y$.
Assume that the distance between two points in the neighborhoods of $x$ and $y$ admits the following second-order Taylor expansion:
$$d(\exp_x(\varepsilon v),\exp_y(\varepsilon w))=d(x,y)\left[\begin{array}{l}1+\varepsilon\left(l^{(1)}_iv^i+l^{(2)}_jw^j\right)\\+\frac{\varepsilon^2}{2}\left(q^{(1)}_{i_1i_2}v^{i_1}v^{i_2}+q^{(2)}_{j_1j_2}w^{j_1}w^{j_2}+2q^{(12)}_{ij}v^iw^j\right)+o\left(\frac{\varepsilon^2}{|\ln(\varepsilon)|}\right)\end{array}\right].$$
Then the coarse Ricci curvature between $x$ and $y$ is:
$$\kappa(x,y)=-l^{(1)}_iF^i(x)-l^{(2)}_jF^j(y)-\frac{q^{(1)}_{i_1i_2}A^{i_1i_2}(x)+q^{(2)}_{j_1j_2}A^{j_1j_2}(y)}{2}+\tr\left(\sqrt{A^{i_1i_2}(x)q^{(12)}_{i_2j_1}A^{j_1j_2}(y)q^{(12)}_{i_3j_2}}\right).$$
\end{thm}
Here the matrix $S^{i_1}{}_{i_3}=A^{i_1i_2}(x)q^{(12)}_{i_2j_1}A^{j_1j_2}(y)q^{(12)}_{i_3j_2}$ is diagonalizable with non-negative eigenvalues, since it is the product of two symmetric non-negative matrices, so $S^{i_1}{}_{i_2}$ admits an unique diagonalizable square root $R^{i_1}{}_{i_2}$ with non-negative eigenvalues, and the last term of the formula is simply $R^i{}_i$.

\begin{rmk} We don't assume here that $d$ is the usual geodesic distance on the manifold $\M$, but only that it admits a nice second order Taylor expansion.
For example, we can take $d$ the Euclidean distance on the sphere $S^n$ embedded in $\Rl^{n+1}$.
\end{rmk}

\noindent\textbf{Proof: }
The idea is to approximate the distributions $P_x^t$ and $P_y^t$ for small $t$ by Gaussian distributions in the tangent spaces $\T_x\M$ and $\T_y\M$, and to approximate the distance by its second order Taylor expansion.
We can describe the process $x(t)$ starting at $x$ in the exponential map by the equation:
$$\der X^i(t)=B^i{}_\alpha(X(t))\der W^\alpha(t)+F'^i(X(t))\der t$$
where $W(t)$ is a Brownian motion in $\Rl^n$, $B^{i_1}{}_{\alpha_1}(0)\delta^{\alpha_1\alpha_2}B^{i_2}{}_{\alpha_2}(0)=A^{i_1i_2}(x)$ and $F'^i(0)=F^i(x)$, and $B$ and $F'$ are continuous (because of the continuity of $A$ and $F$) and defined in a neighborhood of $0$.
Keep in mind that $X(t)$ may not be defined for every $t>0$, but we have $x(t)=\exp_x(X(t))$ when it is.
We will approximate $X^i(t)$ by
$$X^{(0)i}(t)=B^i{}_\alpha(0)W^\alpha(t)+tF^i(x),$$
which has the Gaussian law $\mathcal{N}(tF(x),tA(x))$.
For small $t$, the ball $K_t$ of radius $\sqrt{(2A^{i_1i_2}(x)g_{i_1i_2}(x)+2)t|\ln(t)|}$ of $\T_x\M$ is included in the definition domain of $B$ and $F'$.
We will show that $X(s)$ remains in $K_t$ for $0\leq s\leq t$ with probability $1-o(t)$.
Let $T_t$ be the exit time of $X(s)$ from $K_t$, and
$$X_t(s)=\left\{\begin{array}{ll}X(s)&\text{if  }s\leq T_t\\
X^{(0)}(s)-X^{(0)}(T_t)+X(T_t)&\text{if  }s>T_t.\end{array}\right.$$
We want to prove that $\|X_t|_{[0,t]}\|_\infty=\sup_{s\in[0,t]}\|X_t(s)\|\leq\sqrt{(2A^{i_1i_2}(x)g_{i_1i_2}(x)+2)t|\ln(t)|}$ with probability $1-o(t)$.

We first prove that $\|X_t-X^{(0)}|_{[0,t]}\|_\infty=o(\sqrt{t|\ln(t)|})$ with probability $1-o(t)$.
We have $\der(X_t-X^{(0)})(s)=\mathbf{1}_{s<T_t}[(B^i{}_\alpha(X(s))-B^i{}_\alpha(0))\der W^\alpha(s)+(F'^i(X(s))-F^i(x))\der s]$.
Because of the continuity of $B$ and $F'$, we have $\der(X_t-X(0))(s)=o(1)\der W(s)+o(1)\der s$.
We have $\|\int_0^so(1)\der u|_{0,t]}\|_\infty=o(t)=o(\sqrt{t|\ln(t)|})$.
For each coordinate of the martingale $U^i(s)=\int_0^s\mathbf{1}_{s<T_t}(B^i{}_\alpha(X_t(u))-B^i{}_\alpha(0))\der W^\alpha(u)$, we will apply the Doob inequality to the sub-martingale $\nep^{\lambda U^i}$.
Indeed, we have $\frac{\der}{\der s}\Mn[\nep^{\lambda U^i(s)}]\leq \frac{c_1(t)}{2}\lambda^2\Mn[\nep^{\lambda U^i(s)}]$, with $c_1(t)=o(1)$ depending on the infinite norm on $K_t$ of $(B^{i_1}{}_{\alpha_1}(X)-B^{i_1}{}_{\alpha_1}(0))\delta^{\alpha_1\alpha_2}(B^{i_2}{}_{\alpha_2}(X)-B^{i_2}{}_{\alpha_2}(0))$ so $\Mn[\nep^{\lambda U^i(t)}]\leq \nep^{\frac{\lambda^2c_1(t)t}{2}}$.
So the Doob inequality implies, by taking $c_2(t)=\sqrt{(3c_1(t))t|\ln(t)|}=o(\sqrt{t|\ln(t)|})$, and taking $\lambda=\pm\frac{c_2(t)}{c_1(t)}$, that $\Pb(\sup_{[0,t]}|U^i(s)|\geq c_2(t))\leq 2\nep^{-\frac{c_2(t)^2}{2tc_1(t)}}=2\nep^{-\frac{3|\ln(t)|}{2}}=2t^{\frac{3}{2}}=o(t)$.
We deduce $\|U|_{[0,t]}\|_\infty\leq o(\sqrt{t\ln(t)})$ with probability $1-o(t)$, so we have the same conclusion for $\|X_t-X^{(0)}|_{[0,t]}\|_\infty$.

We have $\|sF^i(0)\|_\infty=O(t)=o(\sqrt{t|\ln(t)|})$, so it remains to prove that $\|B^i{}_\alpha(0) W^\alpha(s)|_{[0,t]}\|_\infty\leq\sqrt{(2A^{i_1i_2}(x)g_{i_1i_2}(x)+1)t|\ln(t)|}$ with probability $1-o(t)$.
We can suppose we are working with an orthonormal basis of eigenvectors of $A^{i_1i_2}(x)g_{i_2i_3}(x)$.
In this case, using the same method as above, we prove $\Pb(\sup_{[0,t]}B^i{}_\alpha(0)W^\alpha(s)\leq\sqrt{(2\lambda_i+\frac{1}{n})}t|\ln(t)|)$ with probability $1-o(t)$, and then if the inequality is true for all $i$, we get the announced result by summing the squares.

Now we set
$$\bar{X}_t(s)=\left\{\begin{array}{ll} X(s)&\text{if }s\leq T_t\\
0&\text{otherwise.}\end{array}\right.$$
We have $\Mn[d(x(t),\exp_x(\bar{X}_t(t)))]=o(t)$.
Indeed, if $X(t)$ does not exit from $K_t$ ($T_t\geq t$), then the distance is $0$.
If $T_t<t$ (what we have shown to occur with probability $o(t)$), we apply the Markov property, and using the local uniform $L^1$-boundedness assumption, the conditional expectation of $d(x(t),x)$ knowing $(X(T_t),T_t)$ is smaller than $M$ for $t$ small enough.
So $\Mn[d(x(t),\exp_x(\bar{X}(t)))]\leq Mo(t)=o(t)$.
So the Wasserstein distance between the distributions of $x(t)$ and $\exp_x(\bar{X}_t(t))$ is $o(t)$.

Of course, we can do the same for the process starting at $y$, and define $Y^{(0)}$, $K'_t$, $Y_t$ and $\bar{Y}_t$.

We denote $\tilde{d}$ the second order Taylor expansion of the distance:
$$\tilde{d}(X,Y)=d(x,y)[1+l^{(1)}_iX^i+l^{(2)}_jY^j+\frac{1}{2}(q^{(1)}_{i_1i_2}X^{i_1}X^{i_2}+q^{(2)}_{j_1j_2}Y^{j_1}Y^{j_2}+2q^{(12)}_{ij}X^iY^j)]$$
The supremum of $\tilde{d}(X,Y)-d(\exp_x(X),\exp_y(Y))$ over $K_t\times K'_t$ is $o(\frac{\sqrt{t|\ln(t)|}^2}{|\ln(\sqrt{t|\ln(t)|})|})=o(t)$.
So the Wasserstein distance between $\exp_x(\bar{X}_t(t))$ and $\exp_y(\bar{Y}_t(t))$ differs of $o(t)$ from the minimum over all couplings of $\Mn[\tilde{d}(\bar{X}_t(t),\bar{Y}_t(t))]$.

It remains to prove that we have a difference of $o(t)$ between the solutions of the minimization problems of $\Mn[\tilde{d}(\bar{X}_t(t),\bar{Y}_t(t))]$ and $\Mn[\tilde{d}(X^{(0)}(t),Y^{(0)}(t))]$.
The expectation and the covariance of $\bar{X}_t(t)-X^{(0)}(t)$ and $\bar{Y}_t(t)-Y^{(0)}(t)$ are $o(t)$.
Indeed, we have $\bar{X}_t(t)-X^{(0)}(t)=(\bar{X}_t(t)-X_t(t))+(X_t(t)-X^{(0)}(t))$.
The expectation and the covariance of $X_t(t)-X^{(0)}(t)$ are $o(t)$ because this quantity is $\int_0^to(1)\der W(s)+o(1)\der s$ (the stochastic integral is a martingale, so its expectation at time $t$ is $0$).
We have $\bar{X}_t(t)-X_t(t)=0$ when $T_t\geq t$, which occurs with probability $1-o(t)$, and the conditional expectation and covariance of $\bar{X}_t(t)-X_t(t)$ knowing $T_t<t$ are $O(\sqrt{t|\ln(t)|})$ and $O(t|\ln(t)|)$, so the expectation and covariance of $\bar{X}_t(t)-X_t(t)$ are $o(t\sqrt{t|\ln(t)|})$ and $o(t^2|\ln(t)|)$ (so $o(t)$ anyway).
So the expectation and the covariance of $\bar{X}_t(t)-X^{(0)}(t)$ are $o(t)$.

The expectation and the covariance of $X^{(0)}(t)$ and $Y^{(0)}(t)$ are $O(t)$, so for any coupling of $(\bar{X}_t(t),X^{(0)}(t))$ with $(\bar{Y}(T),Y^{(0)}(t))$, we have
\begin{multline*}\Mn[\tilde{d}(\bar{X}_t(t),\bar{Y}_t(t))-\tilde{d}(X^{(0)}(t),Y^{(0)}(t))]=\\
d(x,y)\Mn\left[\begin{array}{c}l_i^{(1)}(\bar{X}_t^i(t)-X^{(0)i}(t))+l_j^{(2)}(\bar{Y}_t^j(t)-Y^{(0)j}(t))\\
+\frac{1}{2}q^{(1)}_{i_1i_2}(\bar{X}_t^{i_1}(t)-X^{(0)i_1}(t))(\bar{X}_t^{i_2}(t)+X^{(0)i_2}(t))\\
+\frac{1}{2}q^{(2)}_{j_1j_2}(\bar{Y}_t^{j_1}(t)-Y^{(0)j_1}(t))(\bar{Y}_t^{j_2}(t)+Y^{(0)j_2}(t))\\
+q^{(12)}_{ij}(\bar{X}_t^i(t)-X^{(0)i}(t))(\bar{Y}_t^j(t)+Y^{(0)j}(t))\\
+q^{(12)}_{ij}(\bar{X}_t^i(t)+X^{(0)i}(t))(\bar{Y}_t^j(t)-Y^{(0)j}(t))\end{array}\right]=o(t).\end{multline*}
The last four terms above are $o(t)$ because of the Cauchy Schwarz inequality, which implies that for every family of random vectors $V_t$ and $W_t$ whose covariance matrices satisfy $\Cov(V_t,V_t)=o(t)$ and $\Cov(W_t,W_t)=O(t)$, we have $\Cov(V_t,W_t)=o(t)$.
So we have proved that $W_1(P_x^t,P_y^t)=\inf\Mn[\tilde{d}(X^{(0)}(t),Y^{(0)}(t))]+o(t)$.
The laws of $X^{(0)}(t)$ and $Y^{(0)}(t)$ are $\mathcal{N}(tF(x),tA(x))$ and $\mathcal{N}(tF(y),tA(y))$, so we have
$$\Mn[\tilde{d}(X^{(0)}(t),Y^{(0)}(t))]=d(x,y)\left[\begin{array}{l}1+t(l^{(1)}_iF^i(x)+l^{(2)}_jF^j(y))+\frac{t}{2}(q^{(1)}_{i_1i_2}A^{i_1i_2}(x)+q^{(2)}_{j_1j_2}A^{j_1j_2}(y))\\
+\frac{t^2}{2}\left(q^{(1)}_{i_1i_2}F^{i_1}(x)F^{i_2}(x)+q^{(2)}_{j_1j_2}F^{j_1}(y)F^{j_2}(y)+2q^{(12)}_{ij}F^i(x)F^j(y)\right)\\
+\Mn[q^{(12)}_{ij}(X^{(0)i}(t)-tF^i(x))(Y^{(0)j}(t)-tF^j(y))]\end{array}\right].$$
We only have to minimize the last term, and the minimum is $-t\tr\left(\sqrt{A^{i_1i_2}(x)q^{(12)}_{i_2j_1}A^{j_1j_2}(y)q^{(12)}_{i_3j_2}}\right)$ according to Lemma \ref{qad} below.

So we have proved that
$$W_1(P_x^t,P_y^t)=d(x,y)\left[1-t\left(\begin{array}{l}-l^{(1)}_iF^i(x)-l^{(2)}_jF^j(y)-\frac{1}{2}(q^{(1)}_{i_1i_2}A^{i_1i_2}(x)+q^{(2)}_{j_1j_2}A^{j_1j_2}(y))\\
+\tr\left(\sqrt{A^{i_1i_2}(x)q^{(12)}_{i_2j_1}A^{j_1j_2}(y)q^{(12)}_{i_3j_2}}\right)\end{array}\right)+o(t)\right]$$
which precisely means that
$$\kappa(x,y)=-l^{(1)}_iF^i(x)-l^{(2)}_jF^j(y)-\frac{1}{2}(q^{(1)}_{i_1i_2}A^{i_1i_2}(x)+q^{(2)}_{j_1j_2}A^{j_1j_2}(y))+\tr\left(\sqrt{A^{i_1i_2}(x)q^{(12)}_{i_2j_1}A^{j_1j_2}(y)q^{(12)}_{i_3j_2}}\right). \square$$

\begin{lem}\label{qad}Let $A^{i_1i_2}$ and $B^{j_1j_2}$ be two symmetric non-negative tensors belonging to $E_1\otimes E_1$ and $E_2\otimes E_2$, with $E_1$ and $E_2$ two finite dimensional $\Rl$-vector spaces, not necessarily of the same dimension.
Let $D_{ij}$ be a tensor belonging to $E^*_1\otimes E^*_2$.
Then the minimum of $\Mn[D_{ij}X^iY^j]$ over all couplings between $X$ of law $\mathcal{N}(0,A)$ and $Y$ of law $\mathcal{N}(0,B)$ is
$$-\tr\left(\sqrt{A^{i_1i_2}D_{i_2j_1}B^{j_1j_2}D_{i_3j_2}}\right).$$
\end{lem}

\noindent\textbf{Proof:} The quantity to be minimized only depends on the covariance $C^{ij}=\Mn[X^iY^j]$ between $X$ and $Y$ (this quantity is $C^{ij}D_{ij}$).
So our problem is equivalent to minimizing $C^{ij}D_{ij}$ over the set of all possible $C$ such that there exists a coupling between $X$ and $Y$ such that the covariance between $X$ and $Y$ is $C$.
Since $X$ and $Y$ are Gaussian, $C$ is the covariance of a coupling between $X$ and $Y$ if and only if
$$\left(\begin{array}{cc}A&C\\
C^T&B\end{array}\right)$$
is a symmetric non-negative matrix (because there exists a Gaussian coupling having this covariance).
This condition is equivalent to $\forall (X^*_i,Y^*_j)\in E^*_1\times E^*_2,X^*_{i_1}A^{i_1i_2}X^*_{i_2}+Y^*_{j_1}B^{j_1j_2}Y^*_{j_2}+2X^*_iC^{ij}Y^*_j\geq0$, which is equivalent to $\forall (X^*_i,Y^*_j),|X^*_iC^{ij}Y^*_j|\leq\sqrt{X^*_{i_1}A^{i_1i_2}X^*_{i_2}Y^*_{j_1}B^{j_1j_2}Y^*{j_2}}$.
In particular, this implies $C\in\Img(A)\otimes\Img(B)$ (just take $X^*\in\Ker(A)$ or $Y^*\in\Ker(B)$ and remember $\Img(A^T)=(\Ker(A))^\bot$).

Let $n_1=\rk(A)$ and $n_2=\rk(B)$ be the ranks of $A$ and $B$, using suitable bases of $\Img(A)$ and $\Img(B)$, we find "square roots" $A'^i{}_\alpha$ and $B'^j{}_\beta$ of $A$ and $B$, in the sense that $A'I^{(n_1)}A'^T=A$ and $B'I^{(n_2)}B'^T=B$ ($I_{(n)}$ the scalar product on a canonical $n$-dimensional Euclidean space and $I^{(n)}$ the associated scalar product in the dual of this space).
Then $A'$ and $B'$ admit left inverses $A'^{-1}$ and $B'^{-1}$ (here we don't necessarily have unicity, we just choose two left inverses).
We set $C'=A'^{-1}CB'^{-1T}$.
We have
\begin{align*}\left(\begin{array}{cc}A&C\\
C^T&B\end{array}\right)\geq0&\Leftrightarrow\left(\begin{array}{cc}A'^{-1}&0\\
0&B'^{-1}\end{array}\right).\left(\begin{array}{cc}A&C\\
C^T&B\end{array}\right).\left(\begin{array}{cc}A'^{-1T}&0\\
0&B'^{-1T}\end{array}\right)\geq0\\
&\Leftrightarrow\left(\begin{array}{cc}I^{(n_1)}&C'\\
C'^T&I^{(n_2)}\end{array}\right)\geq0\end{align*}
(because $A'A'^{-1}$ restricted to $\Img(A)$ is the identity, and likewise for $B'B'^{-1}$).
So we have reduced the problem to the case where $E'_1$ and $E'_2$ are Euclidean spaces of dimensions $n_1$ and $n_2$, with $A=I^{(n_1)}$, $B=I^{(n_2)}$ and $D'=A'^TDB'$ instead of $D$.
There exist two nice orthogonal bases so that the matrix of $D'$ in the associated dual bases has the following form:
$$D'=\left(\begin{array}{cc}\mathrm{diag}(\lambda_1,\dots,\lambda_r)&0\\
0&0\end{array}\right)$$
with $\mathrm{diag}(\lambda_1,\dots,\lambda_r)$ the diagonal matrix with coefficients $\lambda_1,\dots,\lambda_r$, $\lambda_k>0$ and $r=\rk(D')=\rk(ADB)$, and furthermore, we have unicity of the coefficients $\lambda_k$.
This result can be proved thanks to the polar decomposition.
We have
$$\left(\begin{array}{cc}I^{(n_1)}&C'\\
C'^T&I^{(n_2)}\end{array}\right)>0\Leftrightarrow \|C'\|_{op}\leq1$$
with $\|C'\|_{op}$ the operator norm of $C'$ associated with the Euclidean norms, hence the coefficients of $C'$ are greater than or equals to $-1$.
The minimum of $C'^{\alpha\beta}D'_{\alpha\beta}$ is then $-\sum_{k=1}^r\lambda_k$, and is attained when the matrix of $C'$ in the nice bases is
$$C'=\left(\begin{array}{cc}-I_r&0\\
0&C''\end{array}\right)$$
with $\|C''\|_{op}\leq1$, and only for those ones $C'$.

The endomorphism $I^{(n_1)}D'I^{(n_2)}D'^T$ has the eigenvalues $\lambda_k^2$ and $0$ with multiplicity $n_1-r$ (it's matrix in the nice basis is $\mathrm{diag}(\lambda^2_1,\ldots,\lambda^2_r,0,\ldots,0)$).
We have $$I^{(n_1)}D'I^{(n_2)}D'^T=I^{(n_1)}A'^TDB'I^{(n_2)}B'^TD^TA'=I^{(n_1)}A'^TDBD^TA'.$$
For any two matrices $M,N$ of size $p\times q$ and $q\times p$, we have for every $n\in\N$, $\tr((MN)^n)=\tr((NM)^n)$, so $MN$ and $NM$ have the same eigenvalues with the same multiplicity, except for the eigenvalue $0$, where the difference of the multiplicities is $|p-q|$.
So the matrix
$$A'I^{(n_1)}A'DBD^T=ADBD^T$$
also has the eigenvalues $\lambda_k^2$ and $0$ with some multiplicity.
The $\lambda_k$ are then the non-zero eigenvalues of $\sqrt{ADBD^T}$.
So the minimum we were looking for is $-\sum_{k=1}^r\lambda_k=-\tr(\sqrt{ADBD^T})$ ($=-\tr(\sqrt{BD^TAD})$ so the symmetry between $A$ and $B$ is respected, which was not straightforward by looking at the formula). $\square$

The two following remarks provide a good understanding of what the set of the solutions of our minimization problem look like.
\begin{rmk}\label{ec} The set of all possible covariances is convex and compact, and the quantity to minimize is linear, so the minimum is attained at an extremal point of this convex set.
Suppose $n_1\geq n_2$, then in the case of an extremal covariance, the coupling between $X$ and $Y$ has the form $Y^j=M^j{}_iX^i$.
Indeed, $C\mapsto A'^{-1}CB'^{-1T}$ restricted to $\Img(A)\otimes\Img(B)$ is linear and bijective, so $C$ is an extremal covariance if and only if $C'$ is an extremal tensor of operator norm smaller than or equals to $1$.
we know that for any tensor $C'$ there exists two orthogonal bases in which the matrix of $C'$ can be written:
$$C'=\left(\begin{array}{c}\mathrm{diag}(\mu_1,\dots,\mu_{n_2})\\
0\end{array}\right)$$
with $\mu_k\geq 0$.
The operator norm of $C'$ is then $\max_{1\leq k\leq n_2}|\mu_k|$.

So $C'$ is an extremal tensor of norm at most $1$ if and only if $\mu_k=1$ for every $k$.

Indeed, if at least one $\mu_k$ is strictly smaller than $1$, $C'$ is a non-trivial convex combination of the tensors whose matrices in the same basis are
$$\left(\begin{array}{c}\mathrm{diag}(\varepsilon_1,\ldots,\varepsilon_{n_2})\\0\end{array}\right)$$
with $\varepsilon_i=\pm1$, and each of these tensors has operator norm $1$.
And conversely, if $\mu_k=1$, then $C'$ is an extremal tensor of norm at most $1$.
Assume that $C'=tC^{(1)}+(1-t)C^{(2)}$ with $t\in ]0,1[$, and $C^{(1)}$ and $C^{(2)}$ have an operator norm smaller than or eqals to $1$.
Then the matrices of $C^{(1)}$ and $C^{(2)}$ have coefficients smaller than or equals to $1$, so their coefficients on the ``diagonal'' must be $1$.
The coefficients outside the ``diagonal'' are $0$ because the sum of the squared coefficients on each row and each column is less than $1$.
So $C^{(1)}=C^{(2)}=C'$.

If $C$ is an extremal covariance, we have $C'^TI_{(n_1)}C'=I^{(n_2)}$ (just do the product of the matrices in the nice bases).
We set then $M=C^TA'^{-1T}I_{(n_1)}A'^{-1}=B'C'^TI_{(n-1)}A'^{-1}$.
The covariance of $Y-MX$ is
$$\begin{array}{r@{}c@{}l}B-MC-C^TM^T+MAM^T&=&B-[B'C'^TI_{(n_1)}A'^{-1}][A'C'B'^T]\\
{}&{}&-[B'C'^TA'^T][A'^{-1T}I_{(n_1)}C'B'^T]\\
{}&{}&+[B'C'^TI_{(n_1)}A'^{-1}][A'I^{(n_1)}A'^T][A'^{-1T}I_{(n_1)}C'B'^T]\\
{}&=&B-B'C'^TI_{(n_1)}C'B'^T=0.\end{array}$$
So $Y=MX$ as previously said.
\end{rmk}

\begin{rmk}\label{sol} For any solution $C$ of our minimization problem, we have
$$\begin{array}{l}CD^TC=A'C'B'^TD^TA'C'B'^T=A'C'D'^TC'B'^T=A'I^{(n_1)}D'I^{(n_2)}B'^T\\
=A'I^{(n_1)}A'^TDB'I^{(n_2)}B'^T=ADB.\end{array}$$
In particular, we have $(CD^T)^2=ADBD^T$ and $(D^TC)^2=D^TADB$.
If we take $C_0$ the solution which corresponds to $C''=0$, $C_0$ is the unique solution with minimal rank (hence the optimal coupling with "the least correlation" between $X$ and $Y$).
We have $\rk(C_0)=\rk(ADB)=\rk(ADBD^T)$, so $\rk(C_0D^T)\leq\rk(ADBD^T)$, and furthermore $\tr(C_0)D^T=-\tr(\sqrt{ADBD^T})$, hence we have $C_0D^T=-\sqrt{ADBD^T}$, and in a similar way $D^TC_0=-\sqrt{D^TADB}$.
Since $C_0D^TC_0=ADB$, we have $\Img(C_0)\supset\Img(ADB)$ and $\Img(C_0^T)\supset\Img(BD^TA)$, and we have in fact equalities because these matrices have the same rank.
As $ADB$, $ADBD^T$ and $D^TADB$ have the same rank, there exist $E$ and $F$ (which will play the role of $D^{-1T}$) such that $ED^TADB=ADB=ADBD^TF$, and then $C_0$ is given by the formula
$$C_0=-\sqrt{ADBD^T}F=-E\sqrt{D^TADB}.$$

For the other solutions, we have
$$C-C_0=A'\left(\begin{array}{cc}0&0\\
0&C''\end{array}\right)B'^T$$
The condition $\|C''\|_{op}\leq1$ is equivalent to the positivity of
$$\left(\begin{array}{cccc}0&0&0&0\\
0&I_{n_1-r}&0&C''\\
0&0&0&0\\
0&C''^T&0&I_{n_2-r}\end{array}\right)$$
which is equivalent to the positivity of
$$\left(\begin{array}{cc}A-C_0B^{-1}C_0^T&C-C_0\\
(C-C_0)^T&B-C_0^TA^{-1}C_0\end{array}\right)$$
where $A^{-1}=A'^{-1T}I_{(n_1)}A'^{-1}$ and $B^{-1}=B'^{-1T}I_{(n_2)}B'^{-1}$.
But we would find the same results for the products $C_0B^{-1}C_0^T$ and $C_0^TA^{-1}C_0$ by taking $A^{-1}$ and $B^{-1}$ such that $AA^{-1}A=A$ and $BB^{-1}B=B$, so this does not depend on the choice of $A'$, $A'^{-1}$, $B'$ or $B'^{-1}$.

We can split $\Img(A)$ as the direct sum of $\Img(ADB)$ and the orthogonal (for the quadratic form induced by $A$ on $\Img(A)$) of this space (which can be written as $\Img(A)\cap\Ker(BD^T)$).
The two matrices $C_0B^{-1}C_0^T$ and $A-C_0B^{-1}C_0^T$ correspond to the decomposition of $A$ on this two subspaces.
The similar remark is valid for the matrices $C_0^TA^{-1}C_0$ and $B-C_0^TA^{-1}C_0$ with respect to the decomposition of $\Img(B)$ as the sum of $\Img(BD^TA)$ and its orthogonal.
An optimal coupling is then any coupling of $X$ and $Y$ satisfying that the covariance between the orthogonal projections (with respect to $A$ and $B$) of $X$ and $Y$ on $\Img(ADB)$ and $\Img(BD^TA)$ is $C_0$.
\end{rmk}

\subsection{The limit of $\kappa(x,y)$ when $x$ and $y$ are close}

Let us look at what the formula given by Theorem \ref{rxy} for $\kappa(x,y)$ becomes when we take $y=\exp_x(\delta u)$, $d$ the usual geodesic distance on Riemannian manifolds and when $\delta$ tends to $0$.
We have the following result that gives the second order Taylor expansion of the geodesic distance on Riemannian manifolds.

\begin{lem}\label{dld} Let $x\in\M$, $(u,v,w)\in(\T_x\M)^3$ such that $g_{ij}u^iu^j=1$, $y=\exp_x(\delta u)$, $w'\in\T_y\M$ obtained from $w$ by parallel transport along the geodesic $t\mapsto\exp_x(\delta tu)$.
Then we have for fixed small enough $\delta$, the following Taylor expansion in $\varepsilon$:
$$d(\exp_x(\varepsilon v),\exp_y(\varepsilon w'))=\delta\left(1+\frac{\varepsilon}{\delta}u^ig_{ij}(w^j-v^j)+\frac{\varepsilon^2}{2\delta^2}\left(r^{(1)}_{ij}v_iv_j+r^{(2)}_{ij}w^iw^j+2r^{(12)}_{ij}v^iw^j+O(\varepsilon^3)\right)\right)$$
with
$$\begin{array}{l}r^{(1)}_{ij}=g_{ij}-g_{ik}u^ku^lg_{lj}-\frac{\delta^2}{3}R_{kilj}u^ku^l+o(\delta^2)\\
r^{(2)}_{ij}=g_{ij}-g_{ik}u^ku^lg_{lj}-\frac{\delta^2}{3}R_{kilj}u^ku^l+o(\delta^2)\\
r^{(12)}_{ij}=-g_{ij}+g_{ik}u^ku^lg_{lj}-\frac{\delta^2}{6}R_{kilj}u^ku^l+o(\delta^2),\end{array}$$
where $R_{kilj}$ is the Riemann tensor of the manifold, and $r^{(1)}_{ij}u^iv^j=r^{(2)}_{ij}u^iv^j=r^{(12)}_{ij}u^iv^j=r^{(12)}_{ij}u^jv^i=0$ (and not only $o(\delta^2)$).
\end{lem}

\noindent\textbf{Proof :} We will take $\delta$ small enough such that $(x,y)$ does not belong to the cut-locus.
Then the Riemannian distance is smooth on a neighborhood of $(x,y)$.

For the term in $\varepsilon$, the well known fact that the sphere of center $x$ and radius $\delta$ is orthogonal at $y$ to the geodesic joining $x$ to $y$ gives us that the part of this term depending on $w$ is proportional to $g_{ij}u^iw^j$.
A similar argument holds for the term in $\varepsilon$ depending on $v$.
Taking $v$ and $w$ proportional to $u$ give the two constants, so we have the term in $\varepsilon$.

For the term in $\varepsilon^2$, we only show that it does only depend on the orthogonal projections of $v$ and $w$ on the orthogonal of $u$, the proof of the behaviour in $\delta$ being based on tedious calculations.
We define $\Sigma_x$ as the image by the exponential map at $x$ of a small ball of the orthogonal of $u$, and $\Sigma_y$ as the image by the exponential map at $y$ of a small ball of the orthogonal of $u'$.
For $\varepsilon$ small enough, the geodesic between $x_1=\exp_x(\varepsilon v)$ and $y_1=\exp_y(\varepsilon w')$ intersects $\Sigma_x$ and $\Sigma_y$ at $x_2=\exp_x(\varepsilon v_1)$ and $y_2=\exp_y(\varepsilon w'_1)$ (we may have to extend the geodesic of $O(\varepsilon)$ beyond $x_1$ and $y_1$).
We have : $d(x_1,y_1)=\overline{d}(x_1,x_2)+d(x_2,y_2)+\overline{d}(y_2,y_1)$ with $\overline{d}(x_1,x_2)=-d(x_1,x_2)$ if we needed to extend the geodesic beyond $x_1$ and $d(x_1,x_2)$ otherwise, and the same for $\overline{d}(y_2,y_1)$.
We also have $v_1=v_2+O(\varepsilon)$ and $w_1=w_2+O(\varepsilon)$, where $v_2=v-\langle{}u,v\rangle{}u$, $w_2=w-\langle{}u,w\rangle{}u$, and the $O(\varepsilon)$ are orthogonal to $u$.
Since in the exponential map, the variation of the metric is of order $2$, we have $d(x_1,x_2)=\varepsilon\|v-v_1\|(1+O(\varepsilon^2))$, and $d(y_2,y_1)=\varepsilon\|w-w_1\|(1+O(\varepsilon^2))$.
So we get $\overline{d}(x_1,x_2)=-\varepsilon\langle{}u,v\rangle+O(\varepsilon^3)$ and $\overline{d}(y_2,y_1)=\varepsilon\langle{}u,w\rangle+O(\varepsilon^3)$, so we find the terms in $\varepsilon$ we expected, and no terms in $\varepsilon^2$.
As $v_1$ and $w_1$ are orthogonal to $u$, we get $d(x_2,y_2)=d(\exp_x(v_2),\exp_y(w'_2))+O(\varepsilon^3)$.
So the $\varepsilon^2$ term does only depend on $v_2$ and $w_2$ as wanted. $\square$

From Theorem \ref{rxy} and Lemma \ref{dld}, we get:
\begin{thm}\label{rxu} Suppose we have a diffusion process on a manifold $\M$ such that $A$ and $F$ are $\mathcal{C}^1$, $\rk(A)=n$ everywhere, locally uniformly $L^1$-bounded.
Then $\kappa(x,\exp_x(\delta u))$ converges to
$$\kappa(x,u)\mathrel{\mathop{:}=}-u^ig_{ij}u^k\nabla_kF^j+\frac{1}{2}R_{kilj}A^{ij}u^ku^l-\frac{1}{4}\overline{u^i\nabla_iA}^{\alpha\beta}\left(\overline{g^{-1}}\otimes\overline{A}+\overline{A}\otimes\overline{g^{-1}}\right)^{-1}_{\alpha\gamma\delta\beta}\overline{u^j\nabla_jA}^{\gamma\delta}$$
when $\delta$ tends to $0$.

Here, for any $M\in\T_x\M\otimes\T_x\M$, we denote by $\overline{M}$ the canonical projection of $M$ to $(\T_x\M/\mathrm{Vect}(u))\otimes(\T_x\M/\mathrm{Vect}(u))$, and the tensor
$$T_{ijkl}=\left(\overline{g^{-1}}\otimes\overline{A}+\overline{A}\otimes\overline{g^{-1}}\right)^{-1}_{ijkl}$$
is uniquely defined by the relationship:
$$T_{ijkl}\left(\overline{g^{-1}}^{jm}\overline{A}^{kn}+\overline{A}^{jm}\overline{g^{-1}}^{kn}\right)=\delta^m_i\delta^n_l.$$
The contraction $T_{ijkl}M^{jk}$ is the unique matrix $N_{il}$ such that $\overline{A}N\overline{g^{-1}}+\overline{g^{-1}}N\overline{A}=M$.
\end{thm}

\begin{rmk} In the special case $A^{ij}=g^{ij}$, we find the usual curvature of the Bakry-Emery theory:
$$\kappa(x,u)=-\langle u,\nabla_uF\rangle+\frac{1}{2}\mathrm{Ric}(u,u).$$
\end{rmk}

\noindent\textbf{Proof: }The hypothesis that $A$ and $F$ are $\mathcal{C}^1$ gives us that the parallel transport of $A(y)$ and $F(y)$ along the geodesic are $A^{ij}(x)+\delta u^k\nabla_kA^{ij}(x)+o(\delta)=A^{ij}(x)+\delta E^{ij}(\delta)$ where $E^{ij}(\delta)$ tends to $u^k\nabla_kA^{ij}(x)$ when $\delta$ tends to $0$, and $F^i(x)+\delta u^k\nabla_kF^i(x)+o(\delta)$.
The application of Theorem \ref{rxy} and Lemma \ref{dld} leads to
$$\begin{array}{r@{\,}l}\kappa(x,\exp_x(\delta u))=&\frac{u^ig_{ij}}{\delta}(F^j-(F^j+\delta u^k\nabla_kF^j+o(\delta)))\\
{}&+\frac{1}{\delta^2}\left[\begin{array}{l}-(A^{ij}+\frac{\delta}{2}E^{ij}(\delta))(g_{ij}-g_{ik}u^ku^lg_{lj}-\frac{\delta^2}{3}R_{kilj}u^ku^l)\\
+\tr(\sqrt{Ar^{(12)}(\delta)(A+\delta E(\delta))r^{(12)T}(\delta)})\end{array}\right]+o(1).\end{array}$$
The difficult point is to understand the behaviour of the square root when $\delta$ tends to $0$.
The quantity under the square root tends to $(A^{ik}(g_{kj}-g_{kl}u^ku^lg_{lj}))^2$, which is of rank $n-1$ (its kernel is $\mathrm{Vect}(u)$).
The square root of matrices is an analytic function in a neigborhood of matrices with positive eigenvalues.
This is why we quotient the space $\T_x\M$ by $\mathrm{Vect}(u)$ (thanks to Lemma \ref{dld}, we know that $r^{(12)}\in u^\bot\otimes u^\bot$).

We need the second-order Taylor expansion of $\tr(\sqrt{M^2+\varepsilon N})$ with $M$ a diagonalizable matrix with positive eigenvalues.
We have $\sqrt{M^2+\varepsilon N}=M+\varepsilon H+\varepsilon^2K+O(\varepsilon^3)$, so we have $HM+MH=N$ and $H^2+KM+MK=0$.
If we work in a diagonalization basis of $M$ (with $\lambda_{(i)}$ the eigenvalues of $M$), we get:
$$H^i{}_j=\frac{N^i{}_j}{\lambda_{(i)}+\lambda_{(j)}}$$
and
$$K^i{}_j=-\frac{1}{\lambda_{(i)}+\lambda_{(j)}}\sum_k\frac{N^i{}_k}{\lambda_{(i)}+\lambda_{(k)}}\frac{N^k{}_j}{\lambda_{(i)}+\lambda_{(k)}}.$$
So we have:
$$\tr(H)=\sum_i\frac{N^i{}_i}{2\lambda_{(i)}}=\frac{1}{2}\tr(M^{-1}N)$$
and
$$\begin{array}{r@{{}={}}l}\tr(K)&-\sum_{i,j}\frac{N^i{}_jN^j{}_i}{2\lambda_{(i)}(\lambda_{(i)}+\lambda_{(j)})^2}\\
{}&-\sum_{i,j}\frac{N^i{}_jN^j{}_i}{4(\lambda_{(i)}+\lambda_{(j)})^2}(\frac{1}{\lambda_{(i)}}+\frac{1}{\lambda_{(j)}})\\
{}&-\sum_{i,j}\frac{N^i{}_jN^j{}_i}{4\lambda_{(i)}\lambda_j(\lambda_{(i)}+\lambda_{(j)})}\\
{}&-\frac{1}{4}\tr(M^{-1}N((I\otimes M+M\otimes I)^{-1}(M^{-1}N)).\end{array}$$

We only have to apply these results with $M=\overline{A}(g-guu^Tg)$, $\varepsilon=\delta$ and $N=\overline{A}(g-guu^Tg)\overline{E(\delta)}(g-guu^Tg)+\frac{\delta}{6}(\overline{A}R(u)\overline{A}(g-guu^Tg)+\overline{A}(g-guu^Tg)\overline{A}R(u))+o(\delta)$, where $R_{ij}(u)=R_{kijl}u^ku^l\in u^\bot\otimes u^\bot$.
We obtain:
$$\begin{array}{r}\tr(\sqrt{Ar^{(12)}(\delta)(A+\delta E(\delta))r^{(12)T}(\delta)})=\tr(\overline{A}(g-guu^Tg))+\frac{\delta}{2}\tr((\overline{E(\delta)}+\frac{\delta}{3}R(u))(g-guu^Tg))\\
-\frac{\delta^2}{4}\tr(\overline{\nabla_uA}(g-guu^Tg)((I\otimes M+M\otimes I)^{-1}(\overline{\nabla_uA}(g-guu^Tg))))+o(\delta^2).\end{array}$$

We have $\tr(\overline{A}R(u))=\tr(AR(u))$, and the last term can be written $-\frac{\delta^2}{4}\tr(\overline{\nabla_uA}((\overline{A}\otimes\overline{g^{-1}}+\overline{g^{-1}}\otimes\overline{A})^{-1}\overline{\nabla_uA}))$ because the inverse of $g-guu^Tg$ (acting on $\T_x\M/\mathrm{Vect}(u)$) is $\overline{g^{-1}}$.
Replacing this expression of the trace of the square root in the expression of $\kappa(x,\exp_x(\delta u))$ cancels the terms of order $\frac{1}{\delta^2}$ and $\frac{1}{\delta}$, and we get the announced result. $\square$

\begin{rmk} The dependency on $u$ of the last term of the formula for the curvature is generally not quadratic (because of the complicated dependency on $u$ of the tensor $(\overline{A}\otimes\overline{g^{-1}}+\overline{g^{-1}}\otimes\overline{A})^{-1}$), but is always non-positive and greater than or equals to the same expression without the bars (which we would have obtained by using the $W_2$ distance instead of the $W_1$ in the definition on $\kappa$, and this expression without the bars depends on $u$ in a quadratic way).
\end{rmk}

\subsection{Construction of the coupling}

Now we will construct a coupling between the paths of the diffusion process thanks to the optimal coupling in the tangent spaces.
In the case when $A$ is invertible everywhere on $\M$, we have $\rk(A(x)q^{(12)}(x,y)A(y))=n-1$.
According to Remarks \ref{ec} and \ref{sol}, we have two extremal covariances $C^+(x,y)$ and $C^-(x,y)$ in the set of the covariances of optimal couplings, given by the formulas:
$$C^+(x,y)=-\sqrt{A(x)q^{(12)}(x,y)A(y)q^{(12)T}(x,y)}p(x,y)+\frac{1}{\sqrt{u^TA(x)^{-1}uu'^TA(y)^{-1}u'}}uu'^T$$
$$C^-(x,y)=-\sqrt{A(x)q^{(12)}(x,y)A(y)q^{(12)T}(x,y)}p(x,y)-\frac{1}{\sqrt{u^TA(x)^{-1}uu'^TA(y)^{-1}u'}}uu'^T$$
with $y=\exp_x(\delta u)$ with $\delta$ small enough, $u'$ the parallel transport of $u$ and $p(x,y)$ any matrix such that $A(x)q^{(12)}(x,y)A(y)q^{(12)T}(x,y)p(x,y)=A(x)q^{(12)}(x,y)A(y)$.
The extremal covariance $C^+(x,y)$ tends to $A(x)$ when $y$ tends to $x$, whereas $C^-(x,y)$ tends to $A(x)-2\frac{uu^T}{u^TA(x)^{-1}u}$ when $u$ stays fixed and $\delta$ tends to $0$, so the coupling with $C^+(x,y)$ generalizes the coupling by parallel transport, whereas the one with $C^-(x,y)$ generalizes the coupling by reflection introduced by Kendall in \cite{ken}.
Here we will use $C^+$ to construct our coupling for Theorem \ref{cent}, because the behaviour of $C^-$ when $\delta$ tends to $0$ is irregular.

So we can construct a coupling between the paths as a diffusion process on $\M\times\M$ (at least in the neighborhood of the diagonal), whose generator is defined by:
$$\begin{array}{r@{}l}L^+(f)(x,y)={}&\frac{1}{2}[A(x)^{ij}\nabla^2_{(11)ij}f(x,y)+A(y)^{ij}\nabla^2_{(22)ij}f(x,y)\\
{}&{}+2C^{+ij}(x,y)\nabla^2_{(12)ij}f(x,y)]+F^i(x)\nabla_{(1)i}f(x,y)+F^i(y)\nabla_{(2)i}f(x,y).\end{array}$$

The coupling above in the case $A=g^{-1}$ is the one of Theorem \ref{cent}.

\noindent\textbf{Proof of Theorem \ref{cent}: }Let us consider the diffusion process of infinitesimal generator $L^+$, which is well defined outside the cut-locus of $\M$.
In the special case of compact Riemannian manifolds, this is true when $d(x,y)$ is strictly smaller than the injectivity radius.
To get the infinitesimal variation of $d(x(t),y(t))$, we have to compute $L^+(f)$ where $f$ has the special form $f(x,y)=\varphi(d(x,y))$ with $\varphi$ regular enough ($\mathcal{C}^2$). We have:
$$\begin{array}{l}\nabla_{(1)i}f(x,y)=\varphi'(d(x,y))\nabla_{(1)i}d(x,y)\\
\nabla_{(2)i}f(x,y)=\varphi'(d(x,y))\nabla_{(2)i}d(x,y)\\
\nabla^2_{(11)ij}f(x,y)=\varphi'(d(x,y))\nabla^2_{(11)ij}d(x,y)+\varphi''(d(x,y))\nabla_{(1)i}d(x,y)\nabla_{(1)j}d(x,y)\\
\nabla^2_{(12)ij}f(x,y)=\varphi'(d(x,y))\nabla^2_{(12)ij}d(x,y)+\varphi''(d(x,y))\nabla_{(1)i}d(x,y)\nabla_{(2)j}d(x,y)\\
\nabla^2_{(22)ij}f(x,y)=\varphi'(d(x,y))\nabla^2_{(22)ij}d(x,y)+\varphi''(d(x,y))\nabla_{(2)i}d(x,y)\nabla_{(2)j}d(x,y)\end{array}$$
with, according to Lemma \ref{dld}:
$$\begin{array}{l}\nabla_{(1)i}d(x,y)=-g_{ij}(x)u^j(x,y)\\
\nabla_{(2)i}d(x,y)=-g_{ij}(y)u^j(y,x)=g_{ij}(y)u'^j(x,y)\\
\nabla^2_{(11)ij}d(x,y)=d(x,y)q^{(1)}_ij(x,y)\\
\nabla^2_{(12)ij}d(x,y)=d(x,y)q^{(12)}_ij(x,y)\\
\nabla^2_{(22)ij}d(x,y)=d(x,y)q^{(2)}_ij(x,y).\end{array}$$
Thus we get:
$$\begin{array}{r@{}l}L^+f(x,y)={}&\frac{1}{2}\varphi''(d(x,y))\left[\begin{array}{l}A^{ij}(x)g_{ik}(x)u^k(x,y)g_{jl}(x)u^l(x,y)+A^{ij}(y)g_{ik}(k)u^k(y,x)g_{jl}(y)u^l(y,x)\\
+2C^{+ij}(x,y)g_{ik}(x)u^k(x,y)g_{jl}(y)u^l(y,x)\end{array}\right]\\
{}&-d(x,y)\varphi'(x,y)\kappa(x,y).\end{array}$$
Since we have $A=g^{-1}$, we get $C^+(x,y)=C_0(x,y)-u(x,y)u^T(y,x)$, with $C_0\in g^{-1}(x)u(x,y)^\bot\otimes g^{-1}(y)u(y,x)^\bot$.
So the term containing $\varphi''(d(x,y))$ is $0$, which means that the variance of $d(x(t),y(t))$ is $o(t)$ when $t$ tends to $0$.
So $\der d(x(t),y(t))=-d(x(t),y(t))\kappa(x(t),y(t))\der t$.
Then by integration of this equality, we get:
$$d(x(t),y(t))=d(x(0),y(0))\nep^{-\int_0^t\kappa(x(s),y(s))\der s}.\square$$

\subsection{The (H) condition and the curvature $\tilde{\kappa}$}

The variance term of this optimal coupling is generally not $0$ in the case when $A\neq g^{-1}$ (nor a multiple of $g^{-1}$).
So we can try to use another coupling, by replacing $C^+(x,y)$ with $\tilde{C}(x,y)$, which is the optimal covariance (for the distance) under the set of covariances which cancel the variance term of $d$ (if this set is non-empty).

We will prove this set is non-empty if and only if the condition
$$(H)\Leftrightarrow\forall u\in\T\M,u^ig_{jk}u^jg_{lm}u^l\nabla_iA^{km}=0$$
is satisfied.

Indeed, the variance term is always nonnegative, so it may vanish if and only if its minimum is $0$.
This is equivalent, according to Lemma \ref{qad}, to
$$\begin{array}{l}2\tr(\sqrt{A(x)g(x)u(x,y)u^T(y,x)g(y)A(y)g(y)u(y,x)u^T(x,y)g(x)})\\
=u^T(x,y)g(x)A(x)g(x)u(x,y)+u^T(y,x)g(y)A(y)g(y)u(y,x),\end{array}$$
which is equivalent to
$$u^T(x,y)g(x)A(x)g(x)u(x,y)=u^T(y,x)g(y)A(y)g(y)u(y,x)$$
(this is the equality case in the inequality between arithmetic and geometric mean).
Differentiating this condition with respect to $y$ along the geodesic starting at $x$ in the direction $u$ gives the condition $(H)$, and of course the converse implication is obtained by integration.

The hypothesis $(H)$ is a very strong hypothesis: for a given metric, the set of the possible $A$ which are nonnegative and satisfy $(H)$ is a convex cone of finite dimension.
Indeed, $H$ is equivalent to: for every geodesic $\gamma(t)$, $A(\gamma(t))(g^{-1}\dot{\gamma}(t))^{\otimes2}$ is constant.
We choose $x\in\M$, and we take a family of vectors $u_{(k)},k=1,\ldots,\frac{n(n+1)}{2}$ such that $\{u_{(k)}^{\otimes2}\}$ is a basis of the symmetric tensors of $\T^2_x\M$. Then we take $x_(k)=\exp_x(\varepsilon u_{(k)})$, with $\varepsilon$ small enough to have $\|\varepsilon u_{(k)}\|<r$, with $r$ the injectivity radius of $\M$.
Then there exists a ball $B$ centered at $x$ such that for every $y\in B$ and every $k$, there exists an unique minimal geodesic joining $y$ and $x_{(k)}$, with velocity $v_{((k)}$ at $y$, and $\{v_{(k)}^{\otimes2}\}$ is a basis of the symmetric tensors of $\T_y\M$.
The knowledge of $A$ at the points $x_{(k)}$ is sufficient to uniquely determine $A$ on the ball $B$.
For any $z\in\M$, we have $x=\exp_z(v)$ for some $v\in\T_z\M$.
We can find a family of vectors $v_{(k)}\in\T_z\M$ in a neighborhood of $v$ such that $\{v_{(k)}^{\otimes2}\}$ is a basis of the symmetric tensors of $\T^2_z\M$, and $\exp_z(v_{(k)})\in B$.
Then the knowledge of $A$ on the points $x_{(k)}$ uniquely determines $A$ on $\M$.

This argument also shows that $A$ is smooth, and the second order Taylor expansion of $A$ in the neighborhood of a single point is sufficient to determine $A$ one the whole manifold.
The condition $(H)$, and the equations obtained by differentiating it twice show that this Taylor expansion must belong to a subspace of dimension $\frac{n(n+1)^2(n+2)}{12}$.

The following examples give the set of the possible $A$ in the cases when $\M$ is an Euclidean space of dimension $n$, the sphere of dimension $n$ or the hyperbolic space of dimension $n$, providing examples where $(H)$ is satisfied without having $A=g^{ij}$.

\begin{xpl} In all three cases mentionned below, $\M$ can be considered as a submanifold of $E=\Rl^{n+1}$ such that the geodesics are the intersection of $\M$ and a two dimensional vector subspace of $E$. Let $(e_1,\ldots,e_{n+1})$ be the canonical basis of $E$ and $(e^*_1,\ldots,e^*_{n+1})$ be the corresponding dual basis
\begin{itemize}
\item We take $\M$ equal to the affine hyperplane of equation $e^*_{n+1}(x)=1$, equipped with the Euclidean metric $\sum_{i=1}^{n}{e^*_i}^2$ in the first case.
\item We put the scalar product $s=\sum_{i=1}^{n+1}{e^*_i}^2$ on $E$, and we take $\M$ equal to the sphere $s(x,x)=1$, equipped with the metric induced by $s$ in the second case
\item We put the quadratic form $q=\sum_{i=1}^n{e^*_i}^2-{e^*_{n+1}}^2$ on $E$, and we take $\M=\{x|q(x,x)=-1\text{ and }e^*_{n+1}(x)>0\}$, equipped with the metric induced by $q$ in the third case.
\end{itemize}
Then we take $T\in E^{*\otimes4}$ any tensor with the same symmetry as a Riemann tensor, that is, $T$ must satisfy $T_{ijkl}=-T_{jikl}=-T_{ijlk}=T_{klij}$ and the Bianchi identity $T_{ijkl}+T_{jkil}+T_{kijl}=0$.
We construct the tensor field $A$ on $\M$ in the following way: let $(x,v)\in \T\M$, we want to have
$$A(x)(g^{-1}v)^{\otimes2}=T_{ijkl}x^iv^jx^kv^l$$
where the sense of the right hand side is given by considering $x$ and $v$ as elements of $E$.
The quadratic dependency in $v$ is trivial, so $A$ is well defined by the previous equation.
Let us consider a unit speed geodesic on $\M$, joining two distinct points $x$ and $y$, and $v$ and $w$ be the speed vectors of the geodesic at points $x$ and $y$.
As said above, the geodesic is included in a two dimensional subspace of $E$, so $(x,v)$ and $(y,w)$ are two bases of this subspace.
Thus there exists a matrix
$$M=\left(\begin{array}{cc}a&b\\c&d\end{array}\right)$$
such that $y=ax+bv$ and $w=cx+dv$.
Then we have $T(y,w,y,w)=\det(M)^2T(x,v,x,v)$ (that is a classical property of the Riemann tensor).
If $l$ is the length of the geodesic, we have:
$$\begin{array}{ll}M=\left(\begin{array}{cc}1&l\\0&1\end{array}\right)&\text{in the case of the Euclidean space,}\\
M=\left(\begin{array}{cc}\cos(l)&\sin(l)\\-\sin(l)&\cos(l)\end{array}\right)&\text{in the case of the sphere,}\\
M=\left(\begin{array}{cc}\ch(l)&\sh(l)\\ \sh(l)&\ch(l)\end{array}\right)&\text{in the case of the hyperbolic space.}\end{array}$$
In each of the three cases, we have $\det(M)=1$.
Thus we have $A(x)(g^{-1}v)^{\otimes2}=A(y)(g^{-1}w)^{\otimes2}$ as wanted.

The linear application $T\mapsto A$ is injective, so the dimension of its image is $\frac{n(n+1)^2(n+2)}{12}$, which is the maximal dimension of the vector space of symmetric tensor fields on $\M$ satisfying the hypothesis $(H)$.
Thus this image is exactly this vector space.
But the tensor fields $A$ which interest us are nonnegative on $\M$, and this implies some restrictions on $T$.
In the cases of the Euclidean space and the sphere, it is true if and only if the ``sectional curvature'' associated to $T$ is nonnegative, whereas in the case of the hyperbolic space, it is true if and only if this ``sectional curvature'' is nonnegative on the planes whose intersection with the cone $q(x,x)=0$ is not $\{0\}$.
\end{xpl}

In the case when $(H)$ is satisfied, the covariances which cancel the variance term of $d$ take the form: $C(x,y)=\tilde{C}_0(x,y)+C'(x,y)$, with
$$\tilde{C}_0(x,y)=-\frac{A(x)g(x)u(x,y)u^T(y,x)g(y)A(y)}{u^T(x,y)g(x)A(x)g(x)u(x,y)}=-\frac{A(x)g(x)u(x,y)u^T(y,x)g(y)A(y)}{u^T(y,x)g(y)A(y)g(y)u(y,x)}$$
and $C'(x,y)$ is such that the big matrix:
$$\left(\begin{array}{cc}A'(x,y)&C'(x,y)\\
C'^T(x,y)&A'(y,x)\end{array}\right)$$
is nonnegative, with $A'(x,y)=A(x)-\frac{A(x)g(x)u(x,y)u^T(x,y)g(x)A(x)}{u^T(x,y)g(x)A(x)g(x)u(x,y)}$.

Using Lemma \ref{qad} again gives us the following expression of $\tilde{\kappa}(x,y)$:
$$\tilde{\kappa}(x,y)=\frac{1}{\delta}\left[\begin{array}{l}-F(y)g(y)u(y,x)-F(x)g(x)u(x,y)-\frac{1}{2}(\tr(A(x)q^{(1)}(x,y))+\tr(A(y)q^{(2)}(x,y)))\\
-\tr(\tilde{C}_0(x,y)q^{(12)T}(x,y))+\tr(\sqrt{A'(x,y)q^{(12)}(x,y)A'(y,x)q^{(12)T}(x,y)})\end{array}\right].$$

We can define $\tilde{\kappa}(x,u)$ as the limit when $\delta$ tends to $0$ of $\tilde{\kappa}(x,\exp_x(\delta u))$. Then we have:
$$\begin{array}{r@{}l}\tilde{\kappa}(x,u)={}&-u^ig_{ij}u^k\nabla_kF^j+\frac{1}{2}A^{ij}R_{ikjl}u^ku^l-\frac{u^ig_{ij}u^k\nabla_kA^{jl}g_{lm}u^n\nabla_nA^{mo}g_{op}u^p}{2u^ig_{ij}A^{jk}g_{kl}u^l}\\
{}&-\frac{1}{4}B^{ij}(A'\otimes(g^{-1}-uu^T)+(g^{-1}-uu^T)\otimes A')^{-1}_{kijl}B^{kl}\end{array}$$
with
$$\begin{array}{l}A'=A-\frac{Aguu^TgA}{u^TgAgu}\\
B^{ij}=\nabla_uA-\frac{\nabla_uAguu^TgA+Aguu^Tg\nabla_uA}{u^TgAgu},\end{array}$$
and as $A'$, $B$ and $g^{-1}-uu^T$ belong to $g^{-1}u^\bot\otimes g^{-1}u^\bot$, we take $(A'\otimes(g^{-1}-uu^T)+(g^{-1}-uu^T)\otimes A')^{-1}$ the unique inverse of $A'\otimes(g^{-1}-uu^T)+(g^{-1}-uu^T)\otimes A'$ in $\left(\T^*_x\M/\mathrm{Vect}(gu)\right)^{\otimes4}$.

And we have the equivalent of Theorem \ref{cent}:
\begin{lem}\label{cent2} If the hypothesis (H) is satisfied, then there exists a coupling between paths such that
$$d(X(t),Y(t))=d(X(0),Y(0))\nep^{-\int_0^t\tilde{\kappa}(X(s),Y(s))\der s}$$
almost surely on the event that for every $0\leq s\leq t$, $d(x',y')^2$ is smooth in a neighborhood of $(X(s),Y(s))$.
\end{lem}

\section{New bounds for the spectral gap}\label{bnd}
The idea to prove Theorems \ref{bharm} and \ref{harm} is to look at the exponential decay of the Lipschitz norm of $P^tf$ when $f$ is Lipschitz with mean $0$ (respect to the reversible probability measure).
Then we use the reversibility assumption to conclude that the variance of $P^tf$ also decreases exponentially fast with the same rate, which is hence a lower bound for the spectral gap.

\noindent\textbf{Proof of Theorem \ref{harm}:} Let $x$ and $y$ be two points of $\M$ such that $d(x,y)<r_i$ where $r_i$ is the injectivity radius of $\M$.
We have $P^tf(y)-P^tf(x)=\Mn[f(Y(t))-f(X(t))]$ for any coupling between paths.
If $f$ is $1$-Lipschitz, then $|f(Y(t))-f(X(t))|\leq d(Y(t),X(t))$, so Lemma $\ref{cent2}$ tells us that $|P^tf(y)-P^tf(x)|\leq d(x,y)\Mn[\nep^{-\int_0^t\tilde{\kappa}(X(s),Y(s))\der s}]$.
For any $\varepsilon>\eta>0$, there exists $\delta>0$ such that for all $(x',y')$ such that $d(x',y')\leq\delta$, we have $\tilde{\kappa}(x',y')\geq\tilde{\kappa}(x')-\eta$, where $\tilde{\kappa}(x')=\inf_{u\in\T_{x'}\M}\tilde{\kappa}(x',u)$.
Taking $T=\frac{\ln(\frac{r_i}{\delta})}{\varepsilon}$, we have $d(X(T),Y(T))\leq\delta$, and then for $t\geq T$ we have $\Mn[\nep^{-\int_0^t\tilde{\kappa}(X(s),Y(s))\der s}]\leq\frac{\delta}{r_i}\Mn[\nep^{-\int_T^{t-T}(\tilde{\kappa}(X(s))-\eta)\der s}]$.
Following what was done in \cite{glw}, we use the Feynman-Kac semigroup $F^t$ generated by $K$, with
$$Kf(x)=\frac{1}{2}A^{ij}(x)\nabla^2_{ij}f(x)+F^i(x)\nabla_i f(x)-(\tilde{\kappa}(x)-\eta)f(x).$$
Indeed we have $\Mn[\nep^{-\int_0^t(\tilde{\kappa}(X(s))-\eta)\der s}]=F^t1(x)$.
The Lipschitz norm of $P^tf$ is at most $\sup_{x\in\M}\frac{\delta}{r_i}\Mn[\nep^{-\int_T^{t-T}\tilde{\kappa}(X(s),Y(s))\der s}]$, for every $t\geq T$.
This quantity is $\sup_{x\in\M}\frac{\delta}{r_i}\Mn_{\delta_x.P^T}F^{t-T}1(y)$, so it is smaller than or equals to $$\frac{\delta}{r_i}\sup_{x\in\M}\|\frac{\der(\delta_x.P^T)}{\der\pi}\|_{L^2(\pi)}\|F^{t-T}1\|_{L^2(\pi)}.$$
Then the Lipschitz norm of $P^tf$ decreases exponentially fast with a better rate than the one of $F^t1$.

The $L^2$-norm of $F^t1$ decreases exponentially with rate $\inf_{h|\int h^2\der\pi=1}\int -hKh\der\pi\geq \inf_{h|\int h^2\der\pi=1}\lambda_1\mathrm{Var}_\pi(h)+\int(\tilde{\kappa}(x)-\eta)h(x)^2\der\pi(x)=\lambda_1+\inf_{h|\int h^2\der\pi=1}\int(\tilde{\kappa}(x)-\eta)h(x)^2\der\pi(x)-\lambda_1(\int h\der\pi)^2$.
The method of Lagrange multiplicators suggests to take $$h(x)=\frac{c}{\tilde{\kappa}(x)-\eta+\alpha},$$ 
with $\alpha$ such that
$$\frac{1}{\lambda_1}=\int\frac{\der\pi(x)}{\tilde{\kappa}(x)-\eta+\alpha}$$
and $c=\frac{1}{\sqrt{\int\frac{\der\pi(x)}{(\tilde{\kappa}(x)-\eta+\alpha)^2}}}$.
With this $h$, we have
$$\int(\tilde{\kappa}(x)-\eta)h(x)^2\der\pi(x)-\lambda_1(\int h\der\pi)^2=-\alpha.$$

This is indeed the minimal $h$ when $\lambda_1$ is at least the harmonic mean $\lambda$ of $\tilde{\kappa}-\inf(\tilde{\kappa})$. We can see it by using Cauchy-Schwarz:
$$\begin{array}{r@{}c@{}l}\int(\tilde{\kappa}(x)-\eta)h(x)^2\der\pi(x)-\lambda_1\left(\int h\der\pi\right)^2&{}\geq{}&\int(\tilde{\kappa}(x)-\eta)h(x)^2\der\pi(x)-\lambda_1\left(\int(\tilde{\kappa}(x)-\eta+\alpha)h^2\der\pi(x)\right)\\
{}&{}&\left(\int\frac{\der\pi(x)}{\tilde{\kappa}(x)-\eta+\alpha}\right)\\
{}&{}={}&-\alpha.\end{array}$$

In the case where $\lambda_1<\lambda$, we take $\alpha=\eta-\inf(\tilde{\kappa})$. This time we get
$$\begin{array}{r@{}c@{}l}\int(\tilde{\kappa}(x)-\eta)h(x)^2\der\pi(x)-\lambda_1\left(\int h\der\pi\right)^2&{}\geq{}&\int(\tilde{\kappa}(x)-\eta)h(x)^2\der\pi(x)-\lambda\left(\int(\tilde{\kappa}(x)-\eta+\alpha)h^2\der\pi(x)\right)\\
{}&{}&\left(\int\frac{\der\pi(x)}{\tilde{\kappa}(x)-\eta+\alpha}\right)+(\lambda-\lambda_1)\int(h\der\pi)^2\\
{}&{}\geq{}&-\alpha.\end{array}$$
A minimizing sequence $h_i(x)$ can be, for example, a sequence such that $h_i^2$ tends to a Dirac at a point where the minimum of $\tilde{\kappa}$ is reached.

In both cases, the exponential decay rate for zero-mean Lipschitz functions is at least $\lambda_1-\alpha$, then by density of the Lipschitz functions on $L^2(\pi)$ and by the reversibility assumption, the exponential decay rate for zero-mean $L^2(\pi)$ functions (which is equal to $\lambda_1$) is also at least $\lambda_1-\alpha$.
Thus $\alpha$ is nonnegative, which means that $\lambda_1$ is at least the harmonic mean of $\tilde{\kappa}-\eta$, so letting $\eta$ tend to $0$ yields the result. $\square$

Purely analytical methods can also be used to prove this result in the case $A^{ij}=g^{ij}$, and they also work when $\inf_{x\in\M}\kappa(x)=0$.
\begin{lem}\label{fornulle}Let $f$ be a regular enough ($\mathcal{C}^3$) function from $\M$ to $\Rl$. Then we have
$$\left.\frac{\der}{\der t}\right|_{t=0}\|\nabla P^tf\|^2=h(2L(h)+u_kg^{kl}\nabla_lA^{ij}\nabla_ihu_j)+h^2\left(\begin{array}{l}u_kg^{kl}\nabla_lA^{ij}\nabla_iu_j-A^{ij}g^{kl}\nabla_iu_k\nabla_ju_l\\+A^{ij}R_{lij\alpha}g^{\alpha\beta}u_\beta g^{kl}u_k+2u_kg^{kl}\nabla_lF^iu_i\end{array}\right)$$
where $h=\|\nabla f\|$ and $\nabla_kf=hu_k$.
\end{lem}

\noindent\textbf{Proof:} We have $\left.\frac{\der}{\der t}\right|_{t=0}\|\nabla P^tf\|^2=2\nabla_kfg^{kl}\nabla_l(Lf)$, and
$$\begin{array}{r@{{}={}}l}\nabla_l(Lf)&\frac{1}{2}(\nabla_lA^{ij}\nabla^2_{ij}f+A^{ij}\nabla^3_{lij}f)+\nabla_lF^i\nabla_if+F^i\nabla^2_{li}f\\
{}&\frac{1}{2}(\nabla_lA^{ij}\nabla^2_{ij}f+A^{ij}(\nabla^3_{ijl}f+R_{lij\alpha}g^{\alpha\beta}\nabla_\beta f))+\nabla_lF^i\nabla_if+F^i\nabla^2_{il}f.\end{array}$$
Differentiating $\nabla_if=h u_i$, we get $\nabla^2_{ij}f=\nabla_ihu_j+h\nabla_iu_j$, and $\nabla^3_{ijl}f=\nabla^2_{ij}hu_l+\nabla_jh\nabla_iu_l+\nabla_ih\nabla_ju_l+h\nabla^2_{ij}u_l$. So we get:
$$\left.\frac{\der}{\der t}\right|_{t=0}\|\nabla P^tf\|^2=hg^{kl}u_k\left[\begin{array}{l}\nabla_lA^{ij}(\nabla_ihu_j+h\nabla_iu_j)\\+A^{ij}\left(\begin{array}{l}\nabla^2_{ij}hu_l+\nabla_jh\nabla_iu_l+\nabla_ih\nabla_ju_l\\+h\nabla^2_{ij}u_l+hR_{lij\alpha}g^{\alpha\beta}u_\beta\end{array}\right)\\+2h\nabla_lF^iu_i+2F^i(\nabla_ihu_l+h\nabla_iu_l)\end{array}\right].$$
Differentiating $g^{kl}u_ku_l=1$ gives $g^{kl}u_k\nabla_ju_l=0$ and $g^{kl}(\nabla_iu_k\nabla_ju_l+u_k\nabla^2_{ij}u_l)=0$, so using these relationships, the above expression can be simplified to get the formula given in Lemma \ref{fornulle}. $\square$

\noindent\textbf{Proof of Theorem \ref{harm3}:} We first prove the Theorem in the case $n'=\infty$ and $c=0$, in which case we get the result of Theorem \ref{harm}. Indeed, in this case, the optimal $\rho(x)$ is nothing but $\kappa(x)=\inf_{u\in\T_x\M}\kappa(x,u)$.

Let $f$ be an eigenfunction of $L$ for the eigenvalue $-\lambda_1$. With the previous notation for $h$ and $u$, we have:
$$\begin{array}{r@{}c@{}l}-2\lambda_1\|h\|^2_{L^2(\pi)}=\left.\frac{\der}{\der t}\right|_{t=0}\|\nabla P^tf\|^2_{L^2(\pi)}&{}={}&\int2h(x)L(h)(x)-2h(x)^2\kappa(x,g^{-1}u(x))\\
{}&{}&-h(x)^2A^{ij}(x)g^{kl}(x)\nabla_iu_k(x)\nabla_ju_l(x)\der\pi(x)\\
{}&{}\leq{}&-2\lambda_1(\int h(x)^2\der\pi(x)-(\int h(x)\der\pi(x))^2)\\
{}&{}&-2\int\kappa(x)h(x)^2\der\pi(x)+0\end{array}$$
where the inequality $\int hL(h)\leq \lambda_1 \mathrm{Var}(h)$ is due to the reversibility assumption. By Cauchy-Schwarz, we have
$$(\int h(x)\der\pi(x))^2\leq \int\frac{\der\pi(x)}{\kappa(x)}\int\kappa(x)h(x)^2\der\pi(x).$$
Finally we get:
$$\int\kappa(x)h^2(x)\der\pi(x)(\lambda_1\int\frac{\der\pi(x)}{\kappa(x)}-1)\geq0,$$
and if $\int\frac{\der\pi(x)}{\kappa(x)}<+\infty$, then $\int\kappa(x)h^2(x)>0$, because $f$ is nonconstant, so $h$ can't be $0$ almost everywhere.
So we have
$$\lambda_1\geq\frac{1}{\int_\M\frac{\der\pi(x)}{\kappa(x)}}.$$

In the general case, we have $n'\geq n$ and the optimal $\rho$ is given by:
$$\rho(x)=\frac{1}{2}\inf_{u\in\T_x\M,\|u\|=1}\mathrm{Ric}(u,u)+\nabla^2_{u,u}\varphi-\frac{(\nabla_u\varphi)^2}{n'-n}.$$
So we have $\rho\leq\kappa$.
Then with the previous notation, we have:
$$\lambda_1(\int_\M h(x)\der\pi(x))^2-\int_\M\rho(x)h^2(x)\der\pi(x)\geq 0$$
because we have just shown the same with $\kappa$ instead of $\rho$.

We also have
$$\begin{array}{r@{\:}l}\int_\M\Gamma_2(f)(x)\der\pi(x)&=\int_\M\frac{1}{2}[L(\frac{h^2}{2}-\langle\nabla f,\nabla(Lf)\rangle]\der\pi=0+\frac{\lambda_1}{2}\int_\M h^2\der\pi\\
{}&\geq\int_\M\rho\Gamma(f)+\frac{1}{n'}L(f)^2\der\pi=\frac{1}{2}\int_\M\rho h^2\der\pi+\frac{\lambda_1}{2n'}\int_\M h^2\der\pi.\end{array}$$

Thus for any $\theta\in[0,1]$ we have:
$$(1-\theta)\lambda_1(\int_\M h\der\pi)^2-\int_\M\left(\rho-\theta\lambda_1\frac{n'-1}{n'}\right)h^2\der\pi\geq0.$$
For $\theta=1$, we have $0\leq\int_\M(\lambda_1\frac{n'-1}{n'}-\rho)h^2\der\pi\leq\lambda_1\frac{n'-1}{n'}-\inf(\rho)\int_\M h^2\der\pi$, this proves the Bakery--\'Emery bound:
$$\lambda_1\geq\frac{n'}{n'-1}\inf(\rho).$$
So for any $c\in[0,\inf(\rho)]$, we take $\theta=\frac{n' c}{(n'-1)\lambda_1}\in[0,1]$.
By Cauchy--Schwarz, we get
$$(1-\theta)\lambda_1\int(\rho-c)h^2\der\pi\int\frac{\der\pi}{\rho-c}-\int(\rho-c)h^2\der\pi\geq0$$
Thus we get $(\lambda_1-c\frac{n'}{n'-1})\int\frac{\der\pi}{\rho-c}-1\geq0$, which leads to the desired result.$\square$

\end{document}